\documentclass[a4paper]{article}
\usepackage{amsmath, amstext, amsxtra, amsthm, amscd, amsgen,amsbsy, amsopn, amsfonts, latexsym, amssymb, amscd}
\usepackage{marvosym}
\usepackage{color}

\begin{document}



\def\wrt{with respect to}

\def\proofend{\hbox to 1em{\hss}\hfill $\blacksquare $\bigskip }

\newtheorem{theorem}{Theorem}[section]
\newtheorem{proposition}[theorem]{Proposition}
\newtheorem{lemma}[theorem]{Lemma}
\newtheorem{remark}[theorem]{Remark}
\newtheorem{remarks}[theorem]{Remarks}
\newtheorem{definition}[theorem]{Definition}
\newtheorem{corollary}[theorem]{Corollary}
\newtheorem{example}[theorem]{Example}
\newtheorem{assumption}[theorem]{Assumption}
\newtheorem{problem}[theorem]{Problem}
\newtheorem{question}[theorem]{Question}
\newtheorem{conjecture}[theorem]{Conjecture}
\newtheorem{rigiditytheorem}[theorem]{Rigidity Theorem}
\newtheorem{mainlemma}[theorem]{Main Lemma}
\newtheorem{claim}[theorem]{Claim}

\def\Z{{\mathbb Z}}
\def\R{{\mathbb R}}
\def\Q{{\mathbb Q}}
\def\C{{\mathbb C}}
\def\N{{\mathbb N}}
\def\H{{\mathbb H}}
\def\Zp #1{{\mathbb Z }/#1{\mathbb Z}}

\def\sec{{\rm sec}}
\def\diam{{\rm diam}}

\def\paperref#1#2#3#4#5#6{\text{#1:} #2, {\em #3} {\bf#4} (#5)#6}
\def\bookref#1#2#3#4#5#6{\text{#1:} {\em #2}, #3 #4 #5#6}
\def\preprintref#1#2#3#4{\text{#1:} #2 #3 (#4)}

\def\diff#1{{\mathrm{Diff}} (#1)}
\def\riemmetric#1#2{{\cal R}_{#2}(#1)}

\makeatletter
\renewcommand\subsection{\@startsection{subsection}{2}{\z@}%
                                     {-3.25ex\@plus -1ex \@minus -.2ex}%
                                     {-1.5ex \@plus .2ex}%
                                     {\normalfont\large\bfseries}}
\renewcommand\subsubsection{\@startsection{subsubsection}{3}{\z@}%
                                     {-3.25ex\@plus -1ex \@minus -.2ex}%
                                     {-1.5ex \@plus .2ex}%
                                     {\normalfont\normalsize\bfseries}}
\makeatother


\author{Anand Dessai\thanks{The first author acknowledges support by SNF grant $200020\_149761$}, Stephan Klaus, and Wilderich Tuschmann}

\title{Nonconnected Moduli Spaces\\ of Nonnegative Sectional Curvature Metrics\\
on Simply Connected Manifolds}

\maketitle


\begin{abstract}
We show that
in each dimension $4n+3$, $n\ge 1$, there exist infinite sequences
of closed smooth simply connected manifolds $M$ of pairwise distinct homotopy type
for which the moduli space of Riemannian metrics with nonnegative sectional curvature
has infinitely many path components.
Closed manifolds with these properties were known before only in dimension seven,
and our result also holds for moduli spaces of Riemannian metrics with positive Ricci curvature.
Moreover, in conjunction with work of Belegradek, Kwasik and Schultz,
we obtain that for each such $M$ the moduli space of complete nonnegative sectional curvature metrics on the open simply connected manifold $M\times\R$ also has
 infinitely many path components.
\end{abstract}

\noindent
\section{Introduction}\label{intro}

A central question in Riemannian geometry concerns the existence of complete metrics on smooth manifolds
which satisfy certain prescribed curvature properties such as, e.g., positivity of scalar or Ricci curvature,
nonnegativity or negativity of sectional curvature, etc.
On the other hand, once the respective existence question is settled, there is an equally important second one,
namely: How `many' metrics of the given type are there, and
how `many' different geometries
of this kind does the manifold actually allow?

To answer these questions, one is naturally led to study the corresponding
spaces of metrics that satisfy the curvature characteristics under investigation,
as well as their respective moduli spaces, i.e., the quotients of these spaces
by the action of the diffeomorphism group given by pulling back metrics.
When equipped with the topology
of smooth convergence on compact subsets and the respective quotient topology, the topological properties and
complexity of these objects constitute the appropiate means to measure the `number'
of different metrics and geometries, and we will adapt to this viewpoint.

\newpage
There has been much activity and profound progress on these issues in the last decades, compare, e.g.,
[BER], [BFK], [BG], [BH], [BHSW],
[BKS], [Ca], [Ch], [CM], [CS], [FO1], [FO2], [FO3], [Ga], [GL], [Hit], [HSS], [LM], [KPT], [KS], [Wa1],
[Wa2], [Wa3], [Wr].
Since the very beginning, special importance has been given to the study
of connectedness properties of spaces and moduli spaces of metrics with lower curvature bounds
on closed as well as open manifolds. As these are also the main issue of the present note,
let us now shortly review and comment on the most relevant developments in this respect,
starting with closed manifolds:

Hitchin ([Hit]) showed first that there are closed manifolds
for which the space of positive scalar curvature metrics is disconnected.
Unfortunately, there is no explicit information on the actual number
of path components available in Hitchin's work. But
for the standard spheres of dimension $4n+3$, where $n\ge 1$,
Carr ([Ca]) then proved that
their respective spaces of positive scalar curvature metrics
indeed have an {\it infinite} number
of path components. Since $\pi _0({\rm Diff}(S^{4n+3}))$ is finite, the corresponding statement also holds for their moduli spaces ([LM, IV \S 7]).

Results about the disconnectedness of moduli spaces
of Riemannian metrics on a wide variety of examples are due to Kreck and Stolz ([KS]).
They showed that for every
closed simply connected spin manifold $M$ of dimension $4n+3$, where $n\ge 1$,
with vanishing
real Pontrjagin classes and $H^1 (M; \Z /2)=0$, the moduli space of metrics of positive scalar curvature, if not empty, has infinitely
many path components.
To achieve this, they introduced in this setting an invariant, the so-called
$s$- (or also, nowadays
Kreck-Stolz) invariant, whose absolute value is constant on path components of the moduli space.

Using the $s$-invariant together with work of Wang and Ziller about
Einstein metrics on principal torus bundles over products of K\"ahler-Einstein manifolds ([WZ]),
Kreck and Stolz went on to show that there are closed seven-manifolds for which the moduli space
of metrics with positive Ricci curvature has infinitely many path components,
and that in this dimension there are also closed manifolds which exhibit a disconnected
moduli space of positive sectional curvature metrics.

Considering different metrics on the Kreck-Stolz examples, it was subsequently
observed in [KPT]
that the Kreck-Stolz result for positive Ricci curvature is actually  also true for moduli
spaces of nonnegative sectional curvature metrics. Notice, however,
that both these results only apply and are indeed confined
to dimension seven.

David Wraith ([Wr]) was then the first to show that there is an infinite number of dimensions with closed manifolds for which the moduli space of Ricci positive metrics has infinitely many path components. Namely, he proved
that for any homotopy sphere $\Sigma^{4n+3}$, $n\ge 1$, bounding a parallelizable manifold,
the moduli space of  metrics
with positive Ricci curvature has an infinite number of path components distinguished by the above Kreck-Stolz invariant.

\

The main result of this note
establishes a corresponding statement
for moduli spaces of
nonnegative sectional curvature metrics
on infinite numbers of manifolds in an infinite range
of dimensions:

\newpage

\begin{theorem}\label{main theorem}
In each dimension $4n+3$, $n\ge 1$, there exist infinite sequences
of closed smooth simply connected manifolds of pairwise distinct homotopy type
for which the moduli space of Riemannian metrics with nonnegative sectional curvature
has infinitely many path components.
\end{theorem}

Moreover, the manifolds figuring here also complement
Wraith's homotopy sphere results. Indeed,
they do actually also show that there is an infinite number
of dimensions in which there exist infinite sequences
of closed smooth simply connected manifolds of pairwise distinct homotopy type
for which the moduli space of  metrics
with positive Ricci curvature has an infinite number of path components. The same holds true if one restricts to  metrics
with positive Ricci and nonnegative sectional curvature.

The manifolds in Theorem \ref{main theorem} are total spaces of principal circle bundles over $\C P^{2n}\times \C P^1$ and agree in dimension $7$ with the examples considered by Kreck and Stolz. Their topology has been studied by Wang and Ziller in their work on Einstein metrics on principal torus bundles mentioned above. These manifolds can also be viewed as quotients of the product of round spheres   $S^{4n+1}\times S^3$ by free isometric circle actions. As in [KPT] we equip these manifolds with the submersion metric. We show that in any fixed dimension $4n+3$ there exist among these Riemannian manifolds infinite sequences with the following properties:
these manifolds belong to a fixed diffeomorphism type and
their metrics have nonnegative sectional and positive Ricci curvature,
but  pairwise distinct absolute $s$-invariants.
Since there are indeed infinite sequences of manifolds of pairwise distinct homotopy types with these properties,
this will imply Theorem~\ref{main theorem}.

\

The study of moduli spaces of nonnegative sectional curvature metrics
 on {\it open} manifolds was inititated in the paper [KPT], with follow-ups in particular in the works
[BFK], [BH] and [BKS].

Using [BKS, Proposition 2.8] one obtains
by means of stabilizing the manifolds in Theorem \ref{main theorem}
with $\R$,
in each dimension $4n$, where $n\ge 2$,
many new examples of simply connected open manifolds
whose moduli spaces of nonnegative sectional curvature metrics have
infinitely many path components. The following result generalizes Corollary 1.11
in [BKS] to an infinite number of dimensions.

\begin{corollary}\label{Corollary}
For each of the manifolds $M$ in Theorem \ref{main theorem},
the moduli space of complete
Riemannian metrics with nonnegative sectional curvature on $M\times\R$
has infinitely many path components.
\end{corollary}

The remaining parts of the present note are structured as follows:
Section $2$ contains the relevant preliminaries needed for the proofs of Theorem \ref{main theorem} and the corollary and the latter are presented in Section \ref{proof}.

\

\noindent
We gratefully acknowledge the support and hospitality of CIRM at Luminy
where a major part of this work was carried out during a research stay within the ``Research in Pairs" program in March 2015. Moreover, it is our pleasure
to thank Igor Belegradek and David Wraith for several helpful comments and the referee for useful suggestions concerning the exposition of the paper.

\newpage

\section{Preliminaries}

\subsection{Spaces and moduli spaces of metrics}\label{moduli spaces}
Let $M$ be an $n$-dimen\-sional smooth  manifold (without boundary), let $S^2 \, T^*M$  denote the second symmetric power of the cotangent bundle of $M$, and let $C^\infty (M,S^2 \, T^*M)$ be the real vector space of smooth symmetric $(0,2)$ tensor fields on $M$. We always topologize $C^\infty (M,S^2 \, T^*M)$ and its subsets with
the smooth topology of uniform convergence on compact
subsets
(compare, e.g., [KrMi], [Ru]). If $M$ is compact, then $C^\infty (M,S^2 \, T^*M)$
is a Fr\'echet space.

The {\it space ${\cal R}(M)$ of all (complete) Riemannian metrics on $M$}
is the subspace  of
$C^\infty (M,S^2 \, T^*M)$ consisting of all sections which are
complete Riemannian metrics on $M$.
Notice that ${\cal R}(M)$ is a convex cone in $C^\infty (M,S^2 \, T^*M)$, i.e., if $a, b > 0$ and $g_1, g_2 \in {\cal R}(M)$, then  $ag_1 + bg_2 \in {\cal R}(M)$.
In particular, ${\cal R}(M)$ is contractible and any open subset of ${\cal R}(M)$ is locally path-connected.

Let us now define moduli spaces of Riemannian metrics.
If $M$ is a finite-dimensional smooth manifold, 
let ${\rm Diff}(M)$ be the group of self-diffeomorphisms of $M$ (which is a Fr\'echet Lie group if $M$ is compact). Then ${\rm Diff}(M)$ acts on ${\cal R}(M)$ by pulling back metrics, i.e., one has the action
$${\rm Diff}(M) \times {\cal R}(M) \to {\cal R}(M), \quad (g,\phi)\mapsto \phi^*(g) \, .$$

The {\it moduli space}
$${\cal M}(M):={\cal R}(M)/\hbox{Diff}(M)$$
of (complete) Riemannian metrics on $M$ is the quotient space of  ${\cal R}(M)$ by the above action of the diffeomorphism group ${\rm Diff}(M)$, equipped with the quotient topology.

Notice that usually ${\rm Diff}(M)$ will not act freely on ${\cal R}(M)$.
Moreover, due to the fact that different Riemannian metrics may have isometry groups of different dimension, the moduli space ${\cal M}(M)$ will in general not have any kind of manifold structure. Note that ${\cal M}(M)$ is locally path-connected since
it is the quotient of the locally path-connected space  ${\cal R}(M)$.
Thus there is no difference between connected and path connected components of open subsets of ${\cal M}(M)$.
Notice also that a lower bound on the number of components in a
moduli space is also a lower bound on the number of components for the respective
space of metrics.

One can similarly form spaces and moduli spaces of metrics satisfying various curvature conditions, as these conditions are invariant under the action of $\hbox{Diff}(M)$.
We will here employ the following notation
(and always tacitly assume that our metrics are complete):
The space of all metrics with positive scalar curvature on $M$
shall be denoted ${\cal R}_{scal >0}(M)$. The corresponding spaces of positive Ricci and nonnegative sectional curvature will be respectively denoted ${\cal R}_{Ric >0}(M)$ and ${\cal R}_{sec \ge 0}(M)$.
The respective moduli spaces will be denoted as ${\cal M}_{scal >0}(M):={\cal R}_{scal >0}(M)/\hbox{Diff}(M)$,
${\cal M}_{Ric >0}(M):={\cal R}_{Ric >0}(M)/\hbox{Diff}(M)$, etc.

Following [KS] we will consider in Section \ref{proof} closed smooth
simply connected manifolds $M$ for which the absolute $s$-invariant is well-defined and constant on path-connected components of the moduli space ${\cal M}_{scal >0}(M)$. Note that in this situation the invariant is also constant on path components of the moduli space ${\cal M}_{Ric >0}(M)$. In other words, two metrics $g_0, g_1\in {\cal R}_{Ric >0}(M)$ with $\vert s\vert (g_0)\neq \vert s\vert (g_1)$ belong to different path components of ${\cal M}_{Ric >0}(M)$.

If, in addition, $g_0, g_1$ are metrics of nonnegative sectional curvature then they also belong to different path components of the moduli space ${\cal M}_{sec\geq 0}(M)$. An elegant way to see this (compare also [BKS, Prop. 2.7])
uses the Ricci flow: Suppose that $\gamma $ is a path
in ${\cal M}_{sec\geq 0}(M)$ with end points represented by $g_0$ and $g_1$. Consider the Ricci flow on ${\cal R}(M)$. Since the Ricci flow is invariant under diffeomorphisms it descends to a local flow on the moduli space. As shown by B\"ohm and Wilking ([BW]),
$\gamma$ evolves under this flow instantly to a path in ${\cal M}_{Ric >0}(M)$. Concatenation of the evolved path and the trajectories of the end points of $\gamma$
then yields a path in ${\cal M}_{Ric >0}(M)$ connecting the end points of $\gamma $, thereby contradicting
$\vert s\vert (g_0)\neq \vert s\vert (g_1)$.

\subsection{The Atiyah-Patodi-Singer Index Theorem}\label{APS-index theorem} An essential ingredient for the definition of the Kreck-Stolz invariant
  is the index theorem of Atiyah, Patodi and Singer for manifolds with boundary. In this section we briefly recall the index theorem for the Dirac and the signature operator.
For more details and the general discussion we refer to [APS1],[APS2].

Let $W$ be a $4m$-dimensional compact spin manifold with boundary $\partial W = M$ and let $g_W$ be a Riemannian metric on $W$ which is a product metric near $M$. Let $g_M$ denote the induced metric on $M$.
Then one can define the Dirac operator $D^+(W,g_W):C^\infty (W,S^+)\to C^\infty (W,S^-)$, where $S^{\pm }$ are the half-spinor bundles. This operator is Fredholm after imposing the APS boundary conditions, i.e. after restricting to spinors $\varphi$ on $W$ for which the restriction to $M$ is in the kernel of the orthogonal projection $P$ onto the space spanned by eigenfunctions for nonnegative eigenvalues [APS1, p. 55]. Its index will be denoted by $\mathrm{ind} \; D^+(W,g_W)$.

In general, the projection $P$ and the APS boundary condition do not depend continuously on the metric due to potential zero eigenvalues of the Dirac operator on the boundary. Hence, the index $\mathrm{ind}\; D^+(W,g_W)\in \Z$ may jump under variations of $g_W$. However, if $g_W(t)$ is a path of metrics as above such that the induced metrics $g_M(t)$ on $M$ have {\em positive} scalar curvature for all $t$ then, by Lichnerowicz's argument, the APS-boundary condition does depend continuously on $t$ and $\mathrm{ind}\; D^+(W,g_W(t))$ is constant in $t$ (see [APSII, p. 417]).

Taking a slightly different point of view,
suppose that $(M,g_M)$ is a $(4m-1)$-dimensional spin manifold of positive scalar curvature which bounds a spin manifold $W$. Let $g_W$ be any extension of $g_M$ to $W$ which is a product metric near the boundary. From the above it follows that $\mathrm{ind}\; D^+(W,g_W)$ is independent of the chosen extension and only depends on the bordism $W$ and the path-component of $g_M$ in the space of metrics on $M$ of positive scalar curvature $\riemmetric {M}{scal>0}$. Moreover the index vanishes if $g_W$ can be chosen to be of positive scalar curvature (see also [KS, Rem. 2.2]).

The index $\mathrm{ind} \; D^+(W,g_W)$ can be computed with the index theorem of Atiyah, Patodi and Singer [APS1, Th. 4.2]. If $(M,g_M)$ has positive scalar curvature one obtains
\begin{equation}\label{APS Dirac}\mathrm{ind}\; D^+(W,g_W)=\int _W {\cal \hat A}(p_1(W,g_W),\ldots ,p_m(W,g_W))-\frac {\eta (D(M,g_M))}2,\end{equation}
where $p_i(W,g_W)$ are the Pontrjagin forms of $(W,g_W)$, $\cal \hat A$ is the multiplicative sequence for the $\hat A$-genus, $D(M,g_M)$ is the Dirac operator on $M$ and $\eta(D(M,g_M))$ is its $\eta$-invariant.

Next we recall the signature theorem for manifolds with boundary. Let $W$ be a $4m$-dimensional compact oriented manifold with boundary $\partial W = M$ and let $g_W$ be a Riemannian metric on $W$ which is a product metric near $M$. Let $g_M$ denote the induced metric on $M$.
Then one can consider the signature operator on $W$ which becomes a Fredholm operator after imposing the APS-boundary conditions. Applying their index theorem to this operator Atiyah, Patodi and Singer proved the following signature theorem [APS1, Th. 4.14]:

\begin{equation}\label{APS signature}\mathrm{sign}\; W = \int _W {\cal L}(p_i(W,g_W))-\eta (B(M,g_M)).\end{equation}
Here ${\cal L}$ is the multiplicative sequence for the $L$-genus and $\eta(B(M,g_M))$ is the $\eta$-invariant of the signature operator on the boundary. Recall that the signature of $W$, $\mathrm{sign}\; W $, is by definition the signature of the (non-degenerate) quadratic form defined by the cup product
on the image of $H^*(W,M)$ in $H^*(W)$. In particular, the right hand side of equation (\ref{APS signature}) is a purely topological invariant.

\subsection{The Kreck-Stolz $s$-invariant}\label{s invariant}

The index of the Dirac operator, the signature and the integrals involving multiplicative sequences considered in the previous subsection depend on the choice of the bordism. However, under favorable  circumstances one can combine the data to define a non-trivial invariant of the boundary itself (cf. [EK], [KS]). In this section we give a brief introduction to the so-called $s$-invariant which was used by Kreck and Stolz in their study of moduli space of metrics of positive scalar curvature.

A starting point in [EK] and [KS] is to consider a certain linear combination ${\cal \hat A} + a_m\cdot \cal L$ of the multiplicative sequences for the Dirac and signature operator. Recall that the degree $\leq 4m$-part of a multiplicative sequence is a polynomial in the variables $p_1,\ldots ,p_m$, where the $p_i$'s are variables which may be thought of as universal rational Pontrjagin classes. By chosing $a_m=1/ (2^{2m+1}\cdot (2^{2m-1}-1))$ in the linear combination above one obtains in degrees $\leq 4m$, a polynomial $N_m(p_1,\ldots ,p_{m-1})$ not involving $p_m$ (cf. also [Hir]).

Now suppose that $W$ is a $4m$-dimensional spin manifold with boundary $M$. In [EK] Eells and Kuiper considered the situation where the real Pontrjagin classes $p_i(W)$, $i<m$, can be lifted uniquely to $H^*(W,M;\R)$ and the natural homomorphism $H^1(W;\Z /2)\to  H^1(M;\Z /2)$ is surjective. This is, for example, the case if $M$ is simply connected with $b_{2m-1}(M)=0$ and $b_{4i-1}(M)=0$, where $0<i<m$.
In this situation Eells and Kuiper showed that the modulo $1$ reduction of
\begin{equation}\label{EK invariant}\langle j^{-1}(N_m(p_1(W),\ldots ,p_{m-1}(W))),[W,M]\rangle- a_m\cdot \mathrm{sign}\; (W)\end{equation}
(or half of it if $m$ is odd) is independent of the choice of $W$ and can be used to detect different smooth structures of $M$. Here $j:H^*(W,M;\R)\to H^*(W;\R )$ is the natural homomorphism induced by inclusion and $\langle \quad ,[W,M]\rangle$ is the Kronecker product with the fundamental homology class.

In [KS] Kreck and Stolz used the APS-index theorem to define an invariant, the so-called $s$-invariant, which refines the Eells-Kuiper invariant. Its relevance for the study of moduli spaces is summarized in Proposition \ref{s moduli proposition} below.

Suppose, as in the last subsection, that $g_W$ is a metric of $W$ which is a product metric near the boundary and let $g_M$ denote the induced metric on $M$. The integral $\int_W N_m(p_1(W,g_W),\ldots ,p_{m-1}(W,g_W))$ can be computed via the APS-index theorem. Suppose $(M,g_M)$ has positive scalar curvature. Then, using equations (\ref{APS Dirac}) and (\ref{APS signature}), one obtains
\small $${\mathrm{ind}\; D^+(W,g_W)+a_m\; \mathrm{sign}\; W =\int _W N_m(p_i(W,g_W))-\frac {\eta (D(M,g_M))}2 -a_m \eta (B(M,g_M)).}
 $$\normalsize
Under favorable  circumstances, e.g., if all real Pontrjagin classes of $M$ vanish, Kreck and Stolz show that the integral $\int _W N_m(p_1(W,g_W),\ldots ,p_{m-1}(W,g_W))$ is equal to
$$\int _M d^{-1}(N_m(p_i(M,g_M))) + \langle j^{-1}(N_m(p_i(W)),[W,M]\rangle,$$ where $d^{-1}(N_m(p_i(M,g_M)))$ is a term which only depends on $(M,g_M)$, i.e., is independent of the choice of the bordism $W$ and the metric $g_W$ (see [KS, p. 829). Collecting the summands involving $(M,g_M)$ one obtains the $s$-invariant of Kreck and Stolz
\begin{equation}\label{s invariant formula} s(M,g_M):= -\frac {\eta (D(M,g_M))}2 -a_m \eta (B(M,g_M)) + \int _M d^{-1}(N_m(p_i(M,g_M)))\end{equation}

Let $t(W):=-(\langle j^{-1}(N_m(p_i(W))),[W,M]\rangle -a_m\cdot \mathrm{sign}\; W)$. Note that the rational number $t(W)$ only depends on the topology of $W$ and reduces modulo $1$ to the negative of the invariant considered by Eells and Kuiper (see equation (\ref{EK invariant})). Now, still asuming that $(M,g_M)$ has positive scalar curvature equations (\ref{APS Dirac}) and (\ref{APS signature}) give
\begin{equation}\label{s equation}s(M,g_M)= \mathrm{ind}\; D^+(W,g_W) +t(W).\end{equation}

Note that the right hand side of (\ref{s equation}) only depends on the bordism $W$ and on the path component of $g_M$ in the space $\riemmetric{M}{scal>0}$ of metrics on $M$ of positive scalar curvature, whereas the left hand side only depends on $(M,g_M)$, i.e. is independent of the chosen bordism.
Hence, both sides depend only on the path component of $g_M$ in $\riemmetric{M}{scal>0}$. If $(W,g_W)$ has positive scalar curvature as well then, by Lichnerowicz's argument, $\mathrm{ind}\; D^+(W,g_W)$ vanishes and $s(M,g_M)$ provides a refinement of the Eells-Kuiper invariant (for all of this see [KS, Section 2]).

If $H^1(M;\Z / 2)=0$, then the spin structure of $M$ is uniquely determined up to isomorphism by the orientation. In this case the absolute $s$-invariant does not change under the action of the diffeomorphism group $\diff{M}$ on $\riemmetric{M}{scal>0}$. In summary, one has

\begin{proposition}[{[KS, Prop. 2.14]}]\label{s moduli proposition}
If $M$ is a closed connected spin manifold of dimension $4m-1$ with vanishing real Pontrjagin classes and $H^1(M;\Z /2)=0$, then $s$ induces a map
$$\vert s\vert :\pi _0(\riemmetric{M}{scal>0}/\diff {M})\to \Q.$$
\end{proposition}

Kreck and Stolz also derive an explicit formula for $s(M,g_M)$ in the case where $M$ is the total space of an $S^1$-principal bundle. We will give this formula in the next section in a particular case.

\subsection{Circle bundles over products of projective spaces}\label{manifolds subsection} This subsection describes the Riemannian manifolds which will be used in the proof of Theorem \ref{main theorem}. As in [KS] we will consider $S^1$-principal bundles over the product of two
complex projective spaces. The total spaces of these bundles have been studied by Wang and Ziller [WZ] in their work about
Einstein metrics on principal torus bundles,
and we will employ some of their cohomological computations.

Let $x$ (resp. $y$) be the positive generators of the integral cohomology ring of $\C P^{2n}$ (resp. $\C P^1$). Let $M_{k,l}$ be the total space of the $S^1$-principal bundle $P$ over $B:=\C P^{2n}\times \C P^1$ with Euler class $c:=lx +ky$, where $k,l$ are coprime positive integers.\footnote{We follow here the notation of Kreck and Stolz. In the notation of Wang and Ziller [WZ] $M_{k,l}$ is denoted by $M^{2n,1}_{l,k}\cong M^{1,2n}_{k,l}$.} We summarize the relevant topological properties of the manifolds $M_{k,l}$ (see also [WZ, Prop. 2.1 and Prop. 2.3]).

\begin{proposition}\label{topology prop}
\begin{enumerate}
\item The spaces $M_{k,l}$ are closed smooth simply connected $(4n+3)$-dimensional manifolds.
\item The integral cohomology ring of $M_{k,l}$ is given by
$$H^*(M_{k,l};\Z )\cong \Z[u,v]/((lv)^{2}, v^{2n+1},uv^2,u^2)$$ where $\deg u=4n+1$
and $\deg v=2$. In particular, $H^*(M_{k,l};\Z )$ does not depend on $k$.
\item The rational cohomology ring of $M_{k,l}$ is isomorphic to the one of $S^{4n+1}\times \C P^1$ and all rational Pontrjagin classes of $M_{k,l}$ vanish.
\item The manifolds $M_{k,l}$ are all formal.
\item For fixed $l$, the manifolds $M_{k,l}$ fall into finitely many diffeomorphism types.
\end{enumerate}
\end{proposition}

\noindent
{\bf Proof:} The first assertion follows from $k$ and $l$ being coprime. The second one can be verified using a spectral sequence argument (cf. [WZ, Prop.~2.1]).
The third statement is a direct consequence of the second.
The fourth statement holds since $H^*(S^{4n+1}\times \C P^1;\Q )$ is intrinsically formal, which follows by direct computation. The last statement is a consequence of
results from Sullivan's surgery theory, compare ([Su], Theorem 13.1 and the proof of Theorem 12.5]).
These imply that if a collection of simply connected closed smooth
manifolds of dimension $\ge 5$ all
have isomorphic integral cohomology rings, the same rational Pontrjagin
classes, and if their minimal model is a formal consequence of their rational
cohomology ring, then there are only finitely many diffeomorphism
types among them. \proofend

Let $L$ denote the complex line bundle associated to the $S^1$-principal bundle $P\to B$ and let $W=D(L)$ be the total space of the disk bundle with boundary $\partial W=S(L)=M_{k,l}$. Note that the tangent bundle of $W$ is isomorphic to $\pi ^*(TB\oplus L)$, where $\pi :W\to B$ is the projection. Hence, the total Stiefel-Whitney class of $W$ is equal to
\begin{equation}\label{SW for W}w(W)\equiv \pi ^*((1+x)^{2n+1}\cdot (1+y)^2\cdot (1+l x +k y))\bmod 2.\end{equation}
Since the restriction of $\pi ^*(L)$ to $M_{k,l}$ is trivial, the tangent bundle of $M_{k,l}$ is stably isomorphic to $\pi ^*(TB)$.

The classes $x$ and $y$ of $H^2(B;\Z )$ pull back under $P\to B$ to $-kv$ and $lv$, respectively (cf. [WZ]). Hence,
 \begin{equation}\label{SW for P}w(M_{k,l})=(1+lv)^2\cdot (1-kv)^{2n+1}.\end{equation}
Restricting formulas (\ref{SW for W}) and (\ref{SW for P}) to $w_2$ and recalling that $M_{k,l}$ and $W$ are simply connected, one obtains the following
\begin{lemma}
$M_{k,l}$ is spin if and only if $k$ is even. In this case $l$ is odd because of $(k,l)=1$,
and, for a fixed orientation, $M_{k,l}$ as well as $W$ admit a unique spin structure.
\proofend
\end{lemma}

From now on we will assume that $k$ is even (and thus $l$ is odd). We equip $W$ with the orientation induced from the standard orientation of $B$ and the complex structure of $L$. From the lemma above we see that $M_{k,l}$ and $W$ admit unique spin structures and that $W$ is a spin bordism for $M_{k,l}$.

Recall from equation (\ref{s equation}) that the computation of the
$s$-invariant of $M_{k,l}$ involves the topological term $t(W):=-(\langle j^{-1}(N_{n+1}(p_i(W)),[W,M]\rangle -a_{n+1}\cdot \mathrm{sign}\; W)$. This term has already been computed by Kreck and Stolz. For the spin bordism $W=D(L)$ as above we thus obtain
\begin{lemma}\label{t(W) lemma}$$t(W)=-\left\langle \frac 1 c \cdot \left( {\cal \hat A}(TB)\cdot \frac {c/2}{\sinh c/2}+a_{n+1}\cdot {\cal L}(TB)\cdot \frac {c}{\tanh c}\right ),[B]\right \rangle $$
\end{lemma}

\noindent
{\bf Proof:} By [KS, Lemma 4.2, part 2]
$$t(W)=-\left\langle {\cal \hat A}(TB)\cdot \frac {1}{2\sinh c/2}+a_{n+1}\cdot {\cal L}(TB)\cdot \frac {1}{\tanh c},[B]\right \rangle +a_{n+1}\cdot \mathrm{sign} (B_c).$$

Here ${\cal \hat A}$ and ${\cal L}$ are the polynomials in the Pontrjagin classes associated to the $\hat A$- and $L$-genus, $a_{n+1}:=1/(2^{2n+3}\cdot (2^{2n+1}-1))$ and $\mathrm{sign} (B_c)$ is the signature of the symmetric bilinear form
\begin{equation}\label{bilinear form}H^{2n}(B)\otimes H^{2n}(B)\to \Q ,\quad (u,v)\mapsto \langle u\cdot v\cdot c, [B]\rangle .\end{equation}
It remains to show that the signature term
$\mathrm{sign} (B_c)$ vanishes. We note that the bilinear form (\ref{bilinear form}) is represented with respect to the basis $(x^n,x^{n-1}y)$ by the matrix $\left (\begin{smallmatrix}
 k & l\\ l & 0\end{smallmatrix}\right )$. Hence, $\mathrm{sign} (B_c)$ vanishes for $l\neq 0$.\proofend

\

\section{Proofs of Theorem \ref{main theorem} and Corollary \ref{Corollary}}\label{proof}
In this section we prove our main theorem and the corollary.
Using the manifolds $M_{k,l}$ from above, we will employ
the absolute Kreck-Stolz invariant $\vert s\vert $ to differentiate between
the path components of the moduli spaces of nonnegative
sectional curvature as well as positive Ricci curvature metrics.

We will first show that each $M_{k,l}$ admits a metric $g_M=g_{k,l}$ of nonnegative sectional curvature and positive Ricci curvature which
is connected in ${\cal M}_{scal >0}(M)$ to a metric that extends to a metric $g_W$ of positive scalar curvature on the associated disk bundle $W$ such that $g_W$ is a product metric near the boundary. This will imply that the $s$-invariant of $(M_{k,l},g_M)$ is equal to the topological term $t(W)$ given in Lemma \ref{t(W) lemma} (see equation (\ref{s equation})).

Next we consider, for fixed $n$, a certain infinite sequence of $S^1$-principal bundles  $P_m\to \C P^{2n}\times \C P^1$ with diffeomorphic total spaces,
where each total space $P_m$ can also be described as the quotient of $S^{4n+1}\times S^3$,
equipped with the standard nonnegatively curved product metric, by the action of an isometric free $S^1$-action.
Notice also that, in fact, there is a free and isometric $T^2=S^1 \times S^1$ - action on $S^{4n+1}\times S^3$
which is
given by the product of the circle Hopf actions on the respective sphere factors, and all
$S^1$-actions that yield the $P_m$ are subactions of this fixed one.

The induced metric $g_m$ on $P_m$ has nonnegative sectional and positive Ricci curvature
by the O'Neill formulas, and an upper sectional curvature bound which is independent of $m$.
We show that the absolute $s$-invariants of the $g_m$ are pairwise different
(for $n=1$ this was already shown by Kreck and Stolz)
and, hence, belong to different connected components of the moduli space of positive scalar (and positive Ricci) curvature metrics.
As explained above,
this implies that the metrics $g_m$ also belong to different connected components of the moduli space of nonnegative sectional curvature metrics (see Subsection \ref{moduli spaces}).

\

Let us now discuss the metrics on $M_{k,l}$ and its associated disk bundle $W$.
\begin{proposition} Each manifold $M:=M_{k,l}$ admits a metric $g_M$ of
simultaneously nonnegative sectional curvature and positive Ricci curvature which
is connected in ${\cal M}_{scal >0}(M)$ to a metric that
extends to a metric $g_W$ of positive scalar curvature on the associated disk bundle $W$ such that $g_W$ is a product metric near the boundary.
\end{proposition}

\noindent
{\bf Proof:} Consider the circle Hopf action on odd dimensional spheres $S^{2r+1}\subset \C^ {r+1}$ given by multiplication of unit complex numbers from the right. This gives rise to an isometric and free action from the right of a two-dimensional torus $T$ on any product of round unit spheres. Here we consider the $T$-action on $S^{4n+1}\times S^3$ with quotient $\C P^{2n}\times \C P^1$. We note that the $T$-orbits are products of geodesics and, hence, totally geodesic flat tori in $S^{4n+1}\times S^3$.

The manifold $M_{k,l}$ can be described as the quotient of $S^{4n+1}\times S^3$ by a certain circle subaction of $T$ (the subaction depends on $(k,l)$). Let $g_M$ be the submersion metric with respect to this action. Then $(M_{k,l}, g_M)$ has nonnegative sectional and positive Ricci curvature by the O'Neill formulas. The (non-effective) action of $T$ on $M_{k,l}$ gives rise also to the $S^1$-principal bundle
$M_{k,l}\to M_{k,l}/T=\C P^{2n}\times \C P^1$. When we equip the base of this fibration with the submersion metric, it is not difficult to see that this Riemannian submersion has totally geodesic fibers.
Now, after shrinking our metric along the fibers if necessary,
following the arguments in [KS, Section 4], we see that the metric can be extended to a metric $g_W$ of positive scalar curvature on the associated disk bundle $W$ such that $g_W$ is a product metric near the boundary.\proofend

Thus, combining the proposition above with Lemma \ref{t(W) lemma}
and equation (\ref{s equation}), it follows that $s(M_{k,l},g_M)$ is equal to
$$-\left\langle \frac 1 c \cdot \left( {\cal \hat A}(TB)\cdot \frac {c/2}{\sinh c/2}+a_{n+1}\cdot {\cal L}(TB)\cdot \frac {c}{\tanh c}\right ),[B]\right \rangle .$$
Let us now compute the invariant $s(k,l):=s(M_{k,l},g_M)$. Recall that $M_{k,l}$ is the total space of the $S^1$-principal bundle over $B=\C P^{2n}\times \C P^1$ with Euler class $c:=lx +k y$ where $x$ and $y$ are the positive generators of the integral cohomology rings of $\C P^{2n}$ and $\C P^1$, respectively, $k,l$ are coprime positive integers, $k$ is even and $l$ is odd.

\begin{proposition} The expression
$s(k,l)/k := s(M^{4n+3}_{k,l},g_M)/k$ is a Laurent polynomial in $l$ of degree $2n$.
\end{proposition}

\noindent
{\bf Proof:} Since the total Pontrjagin class of $B=\C P^{2n}\times \C P^1$
is given by $p(B)=(1+x^2)^{2n+1}$, we see that $s(k,l)$ is equal to
\begin{equation}\label{s polynomial}-\left \langle  \frac 1 c \left(\left(\frac {x/2}{\sinh x/2}\right)^{2n+1}\cdot \frac {c/2}{\sinh c/2}+a_{n+1}\left( \frac {x}{\tanh x}\right )^{2n+1}\cdot \frac {c}{\tanh c}\right ),[B]\right \rangle .\end{equation}
Here
$$\frac {t/2}{\sinh t/2}=1+\sum_{j\geq 1}(-1)^{j}\frac {(2^{2j-1}-1)}{(2j)!\cdot 2^{2j-1}}B_j\cdot t^{2j} =: 1 +\sum _{j\geq 1}\hat a_{2j}\cdot t^{2j}$$ and
$$ \frac {t}{\tanh t}=1+\sum _{j\geq 1}(-1)^{j-1}\frac {2^{2j}}{(2j)!}B_j\cdot t^{2j}=:1 +\sum _{j\geq 1}b_{2j}\cdot t^{2j}$$ are the characteristic power series associated to the $\hat A$- and $L$-genus (cf. [Hir], \S 1.5, [MS], Appendix B). The $B_j$ are the Bernoulli numbers, $B_1=1/6$, $B_2=1/30$, $B_3=1/42$, $B_4=1/30$, $B_5=5/66$,\ldots, and are all non-zero. In particular, $\hat a_{2j}=\frac {1-2^{2j-1}}{2^{4j-1}}\cdot b_{2j}$.

Note that $x^2$ is equal to $c\cdot \frac {lx-ky}{l^2}$ and, hence, divisible by $c$. To evaluate the expression in (\ref{s polynomial}) on the fundamental cycle $[B]$, one first divides out in the inner expression the terms of cohomological degree $4n+4$ by $c$ and then computes the coefficients of $x^{2n}\cdot y$. We claim that the result is $k$ times a Laurent polynomial in $l$ of degree $2n$. To see this we rewrite the formula for $s(k,l)$ in the following form:

$$-s(k,l)=\left \langle \frac 1 c \left(\left(\frac {x/2}{\sinh x/2}\right)^{2n+1}+a_{n+1}\left( \frac {x}{\tanh x}\right )^{2n+1}\right ),[B]\right \rangle $$
$$+\left \langle \left(\frac {x/2}{\sinh x/2}\right)^{2n+1}\cdot \left (\sum _{j\geq 1}\hat a_{2j} c^{2j-1}\right ) +a_{n+1}\left( \frac {x}{\tanh x}\right )^{2n+1}\cdot \left (\sum _{j\geq 1}b_{2j} c^{2j-1}\right ),[B]\right \rangle $$

Using $x^2/c=\frac {lx-ky}{l^2}$ we see that the first summand is a rational multiple of $k\cdot l^{-2}$ and the second summand is of the form $k\cdot p(l)$, where $p(l)$ is a polynomial in $l$ of degree $\leq 2n$. The term of degree $2n$ in $k\cdot p(l)$ is equal to
$$\left \langle \hat a_{2n+2}\cdot c^{2n+1}+a_{n+1}\cdot b_{2n+2}\cdot c^{2n+1},[B]\right \rangle=(2n+1)\cdot k\cdot ( \hat a_{2n+2}+a_{n+1}\cdot b_{2n+2})\cdot l^{2n}$$

Using $a_{n+1}:=1/(2^{2n+3}\cdot (2^{2n+1}-1))$ and $\hat a_{2j}=(1-2^{2j-1})/(2^{4j-1})\cdot b_{2j}$, it then follows that the coefficient of $l^{2n}$ in $p(l)$ is,
indeed, non-zero.\proofend

\noindent
{\bf Proof of Theorem \ref{main theorem}:} We fix a positive odd integer $l_0$ for which the Laurent polynomial $p(l):=s(k,l)/k$ does not vanish.
Note that there are infinitely many of such integers.
Notice also that $s(k,l_0)=k\cdot p(l_0)$ takes pairwise different absolute values for positive even integers $k$. It follows from Sullivan's surgery theory (see Proposition \ref{topology prop}) that there exists an infinite sequence of such integers $(k_m)_m$ with diffeomorphic $M_{k_m,l}$. As before let $P_m:=M_{k_m,l}$ be equipped with the submersion metric.

Recall that $P_m$ has nonnegative sectional and positive Ricci curvature. Since the $P_m$'s have pairwise different absolute $s$-invariant, they belong to pairwise different connected components of the moduli spaces ${\cal M}_{sec \geq 0}(M)$ and ${\cal M}_{Ric >0}(M)$ (see Subsections \ref{moduli spaces} and \ref{s invariant}).

From the classification results of Section $2$ we conclude that there exists, in each relevant dimension, an infinite sequence of manifolds
of pairwise distinct homotopy type which satisfy these properties. \proofend

\noindent
{\bf Proof of Corollary \ref{Corollary}:}
Consider any manifold $M$ as in Theorem \ref{main theorem}.
From the results in Section 2.4 we know that there is
an infinite sequence of nonnegatively curved Riemannian manifolds $(M_k,g_k)$
that are all diffeomorphic to $M$ but have pairwise distinct absolute $s$-invariants. For each $k$ we fix a diffeomorphism $\varphi_k:M\to M_k$ and denote by $\overline g_k$ the pull-back metric $\varphi_k ^*(g_k)$ on $M$. We obtain an infinite sequence of nonnegatively curved metrics $(\overline g_k)_k$ on $M$ which represent pairwise different path components of ${\cal M}_{sec \ge0}(M)$.

Now crossing with $(\R, dt^2)$ gives us nonnegatively curved complete metrics $\overline g_k \times dt^2$ on the simply connected manifold $M\times \R$ with souls $(M, \overline g_k)$ belonging to different path components of ${\cal M}_{sec \ge0}(M)$. By [BKS, Prop. 2.8] the moduli space of complete
Riemannian metrics with nonnegative sectional curvature on $M\times\R$ is homeomorphic to the disjoint union of the moduli spaces of nonnegative sectional curvature metrics of all possible pairwise nondiffeomorphic souls of metrics in ${\cal M}_{sec \ge0}(M\times\R)$.

In particular, the Riemannian manifolds $(M, \overline g_k)\times (\R,dt^2)$ belong to pairwise different path components of the moduli space of complete Riemannian metrics with nonnegative sectional curvature on $M\times\R$.\proofend

\

\subsubsection*{\large References}

\medskip

\begin{description}{\small
\medskip
\item{[APS1]} M. F. Atiyah, V.K. Patodi, I. M. Singer, {\it Spectral asymmetry and Riemannian geometry I}, Math. Proc. Camb. Phil. Soc. {\bf 77} (1975), 43--69.
\medskip
\item{[APS2]} M. F. Atiyah, V.K. Patodi, I. M. Singer, {\it Spectral asymmetry and Riemannian geometry II}, Math. Proc. Camb. Phil. Soc. {\bf 78} (1975), 405--432.
\medskip
\item{[BER]} B. Botvinnik, J. Ebert, O. Randal-Williams, {\it Infinite loop spaces and positive scalar curvature}, Invent. Math. {\bf 209} (2017), 749--835
\medskip
\item{[BFK]} I. Belegradek, F.T. Farrell, V. Kapovitch, {\it Space of non-negatively curved manifolds and pseudoisotopies}, J. Diff. Geom. {\bf 105} (2017), 345--374.
\medskip
\item{[BG]} B. Botvinnik, P. Gilkey, {\it The eta invariant and metrics of positive scalar curvature}, Math. Ann. {\bf 302} (1995), 507--517.
\medskip
\item{[BH]} I. Belegradek and J. Hu, {\it Connectedness properties of the space of complete non-negatively curved planes}, Math. Ann. {\bf 362} (2015), 1273--1286
\medskip
\item{[BHSW]} B. Botvinnik, B. Hanke, T. Schick, M. Walsh, {\it Homotopy groups of the moduli space of metrics of positive scalar curvature}, Geom. Topol. {\bf 14} (2010), 2047--2076.
\medskip
\item{[BKS]} I. Belegradek, S. Kwasik, R. Schultz, {\it Moduli spaces of non-negative sectional curvature and non-unique souls}, J. Diff. Geom. {\bf 89} (2011), 49--86.
\medskip
\item{[BW]} C. B\''ohm and B. Wilking, {\it Nonnegatively curved manifolds with finite fundamental groups admit metrics with positive Ricci curvature}, Geom. Funct. Anal. {\bf 17} (2007), 665--681.
\medskip
\item{[Ca]} R. Carr, {\it Construction of manifolds of positive scalar curvature}, Trans. Amer. Math. Soc. {\bf 307} (1988), 63--74.
\medskip
\item{[Ch]} V. Chernysh, {\it On the homotopy type of the space ${\cal R}^+(M)$}, arxiv:math.GT/0405235.
\medskip
\item{[CM]} F. Cod\'a Marques, {\it Deforming three-manifolds with positive scalar curvature}, Ann. of Math. (2) {\bf 176} (2012), 815--863.
\medskip
\item{[CS]} D. Crowley, T. Schick, {\it The Gromoll filtration, $KO$-characteristic classes and metrics of positive scalar curvature}, Geom. Topol. {\bf 17} (2013), 1773--1790.
\medskip
\item{[EK]} J. Eells, N. Kuiper, {\it An invariant for certain smooth manifolds}, Ann. Math. Pura Appl. {\bf 60} (1962), 93--110.
\medskip
\item{[FO1]} F. T. Farrell, P. Ontaneda,
{\it The Teichm\"uller space of pinched negatively curved metrics on a hyperbolic manifold is not contractible},
Ann. of Math. (2) {\bf 170} (2009), 45--65.
\medskip
\item{[FO2]} F. T. Farrell, P. Ontaneda,
{\it The moduli space of negatively curved metrics of a hyperbolic manifold},
J. Topol. {\bf 3} (2010), 561--577.
\medskip
\item{[FO3]} F. T. Farrell, P. Ontaneda,
{\it On the topology of the space of negatively curved metrics},
J. Diff. Geom. {\bf 86} (2010), 273--301.
\medskip
\item{[Ga]} P. Gajer, {\it Riemannian metrics of positive scalar curvature on compact manifolds with boundary}, Ann. Global Anal. Geom. {\bf 5} (1987), 179--191.
\medskip
\item{[GL]} M. {Gromov} and H.Blaine jun. {Lawson}, {\it Positive scalar curvature and the Dirac operator on complete Riemannian manifolds}, Publ. Math., Inst. Hautes \'Etud. Sci. {\bf 58} (1983), {83--196}.
\medskip
\item{[Hir]} F. Hirzebruch,
{\it Neue topologische Methoden in der algebraischen Geometrie},
Ergebnisse der Mathematik und ihrer Grenzgebiete, Springer, Berlin 1956.
\medskip
\item{[Hit]} N. Hitchin, {\it Harmonic Spinors}, Adv. in Math. {\bf 14} (1974), 1--55.
\medskip
\item{[HSS]} B. Hanke, T. Schick, W. Steimle, {\it The space of metrics of positive scalar curvature}, Publ. Math. Inst. Hautes \'Etudes Sci. {\bf 120} (2014), 335--367.
\medskip
\item{[KPT]} V. Kapovitch, A. Petrunin, W. Tuschmann, {\it Non-negative pinching, moduli spaces and bundles with
infinitely many souls}, J. Diff. Geom. {\bf 71} (2005), 365--383.
\medskip
\item{[KrMi]} A. Kriegl, P. Michor, {\it The convenient setting of global analysis}, Mathematical
Surveys and Monographs {\bf 53}, American Mathematical Society (1997).
\medskip
\item{[KS]} M. Kreck, S. Stolz, {\it Nonconnected moduli spaces of positive sectional curvature metrics}, J. Am. Math. Soc. {\bf 6} (1993), 825--850.
\medskip
\item{[LM]} H.Blaine jun. {Lawson} and M-L {Michelsohn},
{\it Spin geometry}, Princeton University Press 1989.
\medskip
\item{[MS]} J. Milnor and J. Stasheff, {\it Characteristic classes}, Annals of Mathematics Studies no. 76, Princeton University Press 1974.
\medskip
\item{[Ru]} W. Rudin, {\it Functional Analysis}, McGraw-Hill 1973.
\medskip
\item{[Su]} D. Sullivan, {\it Infinitesimal computations in topology},
Publ. Math. IHES {\bf 47} (1977), 269--331.
\medskip
\item{[Wa1]} M. Walsh, {\it Metrics of positive scalar curvature and generalized Morse functions, Part I}, Mem. Amer. Math. Soc. {\bf 209} (2011).
\medskip
\item{[Wa2]} M. Walsh, {\it Cobordism invariance of the homotopy type of the space of positive scalar curvature metrics}, Proc. Amer. Math. Soc. {\bf 141} (2013), 2475--2484.
\medskip
\item{[Wa3]} M. Walsh, {\it H-spaces, loop spaces and the space of positive scalar curvature metrics on the sphere}, Geom. Topol. {\bf 18} (2014), 2189--2243.
\medskip
\item{[Wr]} D. J. Wraith, {\it On the moduli space of positive Ricci curvature metrics on homotopy spheres}, Geom. Topol. {\bf 15} (2011), 1983--2015.
\medskip
\item{[WZ]} M. Wang, W. Ziller, {\it Einstein metrics on principal torus bundles},
J. Differential Geom. {\bf 31} (1990), 215--248.
\medskip
}
\end{description}

\

\

\noindent
{Email\thinspace :\quad}
{klaus@mfo.de, \thinspace anand.dessai@unifr.ch, \thinspace tuschmann@kit.edu}

\

\noindent Mathematisches Forschungsinstitut Oberwolfach (MFO), Schwarzwaldstr. 9-11,
\\ D-77709 Oberwolfach-Walke, Germany

\

\noindent D\'epartement de Math\'ematiques, Chemin du Mus\'ee 23, Facult\'e des sciences,
Universit\'e de Fribourg, P\'erolles, CH-1700 Fribourg, Switzerland

\

\noindent Karlsruher Institut f\"ur Technologie (KIT), Fakult\"at f\"ur Mathematik,
Institut f\"ur Algebra und Geometrie, Arbeitsgruppe Differentialgeometrie,
Englerstr. 2,
\\ D-76131 Karlsruhe, Germany

\

\end{document}